# Alpha-reliable combined mean traffic equilibrium model with stochastic travel times

ZHANG Wen-yi[1], GUAN Wei[1], SONG Li-ying[1], SUN Hui-jun[1]

1. MOE Key Laboratory for Urban Transportation Complex Systems Theory and Technology, Beijing Jiaotong University, Beijing 100044, China

**Abstract:** Based on the reliability budget and percentile travel time (PTT) concept, a new travel time index named combined mean travel time (CMTT) under stochastic traffic network was proposed. CMTT here was defined as the convex combination of the conditional expectations of PTT-below and PTT-excess travel times. The former was designed as a risk-optimistic travel time index, and the latter was a risk-pessimistic one. Hence, CMTT was able to describe various routing risk-attitudes. The central idea of CMTT was comprehensively illustrated and the difference among the existing travel time indices was analysed. The Wardropian combined mean traffic equilibrium (CMTE) model was formulated as a variational inequality and solved via an alternating direction algorithm nesting extra-gradient projection process. Some mathematical properties of CMTT and CMTE model were rigorously proved. In the end, a numerical example was performed to characterize the CMTE network. It is founded that that risk-pessimism is of more benefit to a modest (or low) congestion and risk network, however, it changes to be risk-optimism for a high congestion and risk network.

## 1 Introduction

The reliability-based methodology has been widely used to investigate the traveller’s routing behaviours in stochastic traffic network. Within this framework, a traveller aims not only to minimize the travel time but also to improve the punctual arrival probability (or arrival reliability). Frank [1] employed the reliability-based percentile travel time (PTT) index to measure the quality of a route, and the shortest routes possessed the least PTT. Fan et al. [2] applied PTT into an adaptive routing problem. Note that, when travellers were insatiable, the PTT-minimum routes corresponded to the first-order non-dominated ones according to the stochastic dominance theory [3]. Inspired from this, Nie and Wu [4] explored the PTT-shortest routing problem. Later, the correlations between adjacent links were incorporated [5]. Recently, Nie [6] proposed a percentile user equilibrium model for traffic network analysis.

Another measure to model travel time reliability is the effective travel time (ETT) index. Different from PTT, ETT argues that a traveller reserves a positive safety margin to hedge against the travel time variation in daily trips [7]. Mathematically, ETT comprises the mean travel time (MTT) index and a safety margin which is the product of standard deviation and a scalar called punctuality factor. When the travel times distribute normally, the punctuality factor has a one-to-one correspondence with reliability. However, such relationship does not hold in other cases [8]. Lo and Tung [9] proposed a probabilistic user equilibrium model, where the travellers merely selected the routes experiencing the minimum MTT and satisfying the

**Foundation item:** Project (2012CB725403-5) supported by National Basic Research Program of China; Project (71131001-2) supported by National Natural Science Foundation of China; Project (201170) supported by FANEDD
**Received date:** 2012-07-14; **Accepted date:** 2012-11-29
**Corresponding author:** GUAN Wei, Professor, PhD; Tel: +86-10-51687142; E-mail: weig@bjtu.edu.cn

reliability restrictions. Lo et al. [10] extended the probabilistic user equilibrium model to the travel time budget (TTB) model with degradable link capacities, and the central limit theorem guaranteed the normal distributions of route travel times. When travel times distributed normally TTB was equivalent with PTT, value at risk (VaR) [11] applied widely in finance area, and alpha-reliable travel time index [12]. Based on TTB, many further studies, considering endogenous [13−14] or exogenous random sources [15−16], were explored. Chen and Zhou [17] argued the importance of unreliable aspect in a traveller's routing process. They introduced a mean-excess travel time (METT) index which was defined as the conditional expectation of the TTB-excess travel times. METT was analogous to the conational value at risk (CVaR) [18] used in financial engineering. Compared with TTB, METT actually represented a kind of more pessimistic routing behaviour. Fosgerau and Karlström [19] analysed the value of reliability, measured by the derivative of expected utility with respect to standard deviation, in scheduling activities. They verified that the maximal expected utility was the linear combination of mean and standard deviation for all travel time distributions. Later, Fosgerau and Engelson [20] further reported their study on the value of travel time variance.

Due to the randomness of travel time, someone would like to select the robust routes to avoid the worst possibility. This kind of routing behaviour can be described by the dominant-strategy-based self-routing games. At the equilibrium states, travellers only distribute on the routes with the minimum upper bounds of travel time. In order to obtain these upper bounds, two kinds of methods were developed, i.e., the robust optimization [21] and distribution-free-based method [22].

With the development of ITS technology, more and more (link/route) travel time samples are available. However, the TTB-based methodology is initially developed to describe those risk-pessimistic routing behaviours. METT is a more risk-pessimistic travel time index since METT individuals care merely about the unreliable events with travel times exceeding the budgeted travel time. Analogous pessimistic behaviours may be common in financial investments where investors have to be enough supersensitive and scrupulous to avoid potential investment-ruins because of frequent and huge transactions. In daily trips, travellers also regard the performance of reliable domain, with travel time less than budget, as an important aspect. In order to capture this kind of routing behaviour, a combined mean travel time (CMTT) index is introduced and explored comprehensively in this paper. The remainder of this paper is organized as follows. Section 2 elaborates on the central idea of CMTT. The differences between several existing travel time indices are also analysed. Section 3 presents the corresponding variational inequality model for the combined mean traffic equilibrium (CMTE) and its solution algorithm. In Section 4, a numerical example is conducted to demonstrate the CMTE model and algorithm, as well as to explore the effects of routing risk-attitudes on the network under different demand and capacity degradation levels. Section 5 concludes the study and suggests some future works.

## 2 Combined mean travel time

In this section, to improve the readability, the main abbreviations and notations applied across this study are given in prior. Next, the concept of CMTT is elaborated, and the differences between several travel time indices will be analysed.

### 2.1 Abbreviations

| | |
|---|---|
| VaR | value at risk |
| CVaR | conditional value at risk |
| PTT | percentile travel time |
| MTT | mean travel time |
| ETT | effective travel time |
| TTB | travel time budget |
| METT | mean-excess travel time |
| MBTT | mean-below travel time |
| CMTT | combined mean travel time |
| METE | mean-excess traffic equilibrium |
| MBTE | mean-below traffic equilibrium |
| CMTE | combined mean traffic equilibrium |
| VIP | variational inequality problem |
| NCP | nonlinear complemental problem |
| MAD | modified alternative direction |
| ANTT | average network travel time |
| SA | sensitivity analysis |

## 2.2 Notations

Considering a strongly connected transportation network $G=(N,A)$, where $N$ and $A$ $(a \in A)$ denote the sets of nodes and links, respectively. Let $R$ and $S$ denote a subset of $N$ for which travel demand $q^{rs}$ is generated from origin $r \in R$ to destination $s \in S$. The main notations are listed below, and the unlisted ones will be explained when they are firstly encountered.

| | |
|---|---|
| $k$ | route index |
| $K^{rs}$ | set of routes from $r$ to $s$ |
| $f_k^{rs}$ | flow on route $k$ from $r$ to $s$, $\mathbf{f}=\left(\ldots,f_k^{rs},\ldots\right)^T$ |
| $v_a$ | flow on link $a$ |
| $\alpha$ | confidence level or budgeted reliability |
| $T_k^{rs}$ | random travel time on route $k$ from $r$ to $s$ |
| $t_k^{rs}$ | sample value with respect to $T_k^{rs}$ |
| $T_a$ | random travel time on link $a$ |
| $\xi$ | artificial constant larger than $T_k^{rs}$ |
| $\lambda$ | combined weight of MBTT to estimate CMTT |
| $\xi_k^{rs}(\alpha)$ | alpha-reliable PTT/TTB on route $k$ from $r$ to $s$ |
| $\eta_k^{rs}(\alpha)$ | alpha-reliable METT on route $k$ from $r$ to $s$ |
| $\varphi_k^{rs}(\alpha)$ | alpha-reliable MBTT on route $k$ from $r$ to $s$, $\boldsymbol{\varphi}(\mathbf{f})=\left(\ldots,\varphi_k^{rs}(\alpha),\ldots\right)^T$ |
| $\psi_k^{rs}(\alpha;\lambda)$ | alpha-reliable CMTT on route $k$ from $r$ to $s$, $\boldsymbol{\psi}(\mathbf{f})=\left(\ldots,\psi_k^{rs}(\alpha;\lambda),\ldots\right)^T$ |
| $\pi_\alpha^{rs}$ | minimum alpha-reliable CMTT for the routes in $K^{rs}$, $\boldsymbol{\pi}=(\ldots,\pi_\alpha^{rs},\ldots)^T$ |
| $\mu_k^{rs}$ | expectation/mean of $T_k^{rs}$ |
| $\sigma_k^{rs}$ | standard deviation of $T_k^{rs}$ |
| $\delta_{ak}^{rs}$ | 0-1 indicator variable, and $\delta_{ak}^{rs}=1$ if route $k$ uses link $a$, and 0 otherwise |
| $f(\cdot)$ | operator of probability density function |
| $E(\cdot)$ | expectation operator of random variable |
| $Var(\cdot)$ | variance operator of random variable |
| $\Pr(\cdot)$ | operator of cumulative density function |
| $N(\mu,\sigma)$ | normal distribution indicator with mean $\mu$ and standard deviation $\sigma$ |
| $\Phi(\cdot)$ | cumulative density function for the standard normal distribution |
| $\phi(\cdot)$ | probability density function for the standard normal distribution |
| $\Phi^{-1}(\cdot)$ | inverse of $\Phi(\cdot)$ |

### 2.3 Definition of mean-below travel time

Due to the huge gap in the transaction amount and frequency between everyday travels and normal investments, the travellers cannot behave the same sensitivity towards unreliability as the financial investors do. Therefore, the reliability indices (e.g. CVaR and VaR) used in finance area maybe not well suitable to travel decision-makings. In daily life, one may occasionally suffer from very long duration on an often-used route, yet next time he/she still chooses it. A possible reason is that he/she habitually believes analogous events will emerge on the other routes; it may also attribute to his/her optimistic risk-attitude towards these habitual routes. In other words, he/she believes that the occasional events would not happen to him/her in the limited daily trips. Such kind of overconfidence behaviour was also reported in the studies on the cognitive psychology [23], i.e., people tended to overestimate their own subjective cognitive capacities. Subsequently, a new travel time index named mean-below travel time (MBTT) will be introduced to capture such a kind of risk-optimistic routing behaviour.

**Definition 1.** The alpha-reliable MBTT of a route is the conditional expectation of the corresponding PTT-below travel times, i.e.,

$$\varphi_k^{rs}(\alpha) = E\left(T_k^{rs} : T_k^{rs} \le \xi_k^{rs}(\alpha)\right) \ \forall k \in K^{rs}, r \in R, s \in S \tag{1}$$

Here, $\xi_k^{rs}(\alpha) = \min\left\{\xi : \Pr(T_k^{rs} \le \xi) \ge \alpha\right\}$ can be further rewritten as $\xi_k^{rs}(\alpha) = \left\{\xi : \Pr(T_k^{rs} \le \xi) = \alpha\right\}$. Note that $\xi_k^{rs}(\alpha)$ is equivalent to the alpha-reliable TTB when $T_k^{rs}$ distributes normally. So in the remaining context TTB and PTT will be used interchangeably. Accordingly, Eq. (1) can be formulated as

$$\varphi_k^{rs}(\alpha) = \int_{-\infty}^{\xi_k^{rs}(\alpha)} \frac{t_k^{rs} f(t_k^{rs})}{\Pr(T_k^{rs} \le \xi_k^{rs}(\alpha))} dt_k^{rs} \ \forall k \in K^{rs}, r \in R, s \in S \tag{2}$$

MBTT is regarded as a risk-optimistic travel time index since the MBTT-conducted individuals merely care about the performance of PTT-below reliable domain and believe the PTT-excess unreliable cases will not happen to them in limited trips.

Define $\eta_k^{rs}(\alpha) = E\left(T_k^{rs} : T_k^{rs} \ge \xi_k^{rs}(\alpha)\right)$, which is the alpha-reliable METT explored by Chen and Zhou [17], the following proposition can be obtained.

**Proposition 1.** Equation $\varphi_k^{rs}(\alpha) = \frac{1}{\alpha} E(T_k^{rs}) - \frac{1-\alpha}{\alpha} \eta_k^{rs}(\alpha)$ holds in general for all routes.

**Proof.** Since $E(T_k^{rs}) = \int_{-\infty}^{+\infty} t_k^{rs} f(t_k^{rs}) dt_k^{rs} = \int_{-\infty}^{\xi_k^{rs}(\alpha)} t_k^{rs} f(t_k^{rs}) dt_k^{rs} + \int_{\xi_k^{rs}(\alpha)}^{+\infty} t_k^{rs} f(t_k^{rs}) dt_k^{rs}$

$$= \alpha \int_{-\infty}^{\xi_k^{rs}(\alpha)} \frac{t_k^{rs} f(t_k^{rs})}{\Pr(T_k^{rs} \le \xi_k^{rs}(\alpha))} dt_k^{rs} + (1-\alpha) \int_{\xi_k^{rs}(\alpha)}^{+\infty} \frac{t_k^{rs} f(t_k^{rs})}{\Pr(T_k^{rs} \ge \xi_k^{rs}(\alpha))} dt_k^{rs}$$

$$= \alpha\varphi_k^{rs}(\alpha) + (1-\alpha)\eta_k^{rs}(\alpha)$$

$$\Rightarrow \varphi_k^{rs}(\alpha) = \frac{1}{\alpha}E(T_k^{rs}) - \frac{1-\alpha}{\alpha}\eta_k^{rs}(\alpha)\ \forall k \in K^{rs}, r \in R, s \in S\ \square$$

According to Proposition 1, $\varphi_k^{rs}(1) = E(T_k^{rs})$ is derived, i.e., MBTT changes to be MTT when $\alpha = 1$.

### 2.4 Estimation of mean-below travel time

In order to estimate the alpha-reliable MBTT, an essential step is to obtain the expectation and variation of route travel time. To focus on the essence of MBTT, the link travel times $T_a (a \in A)$ are assumed to be mutually independent. This assumption is also adopted in the other studies (see [15−16]). Accordingly, $\mu_k^{rs} = \sum_{a\in A} E(T_a)\delta_{ak}^{rs}$ and $\sigma_k^{rs} = sqrt\left[\sum_{a\in A} Var\left(T_a\right)\delta_{ak}^{rs}\right]$ for all routes are promised. In addition, assume all the link travel times are bounded. Regardless of the probability distributions of link travel times, applying the Lyapunov condition [24] for the central limit theorem, $T_k^{rs} \sim N\left(\mu_k^{rs}, \sigma_k^{rs}\right)$ can be concluded. Accordingly, the alpha-reliable MBTT, defined by Eq. (1), can be estimated as

$$\varphi_k^{rs}(\alpha) = \mu_k^{rs} - \frac{\sigma_k^{rs}}{\sqrt{2\pi}\alpha}\exp\left(-\frac{\left[\Phi^{-1}(\alpha)\right]^2}{2}\right)\ \forall k \in K^{rs}, r \in R, s \in S \tag{3}$$

The detailed deduction is offered in Appendix A. In Eq. (3), the variability of travel time may result from both exogenous and endogenous disturbances, such as stochastic supply, stochastic demand, etc. The negative scalar before the standard deviation implies MBTT is universally smaller than MTT, making MBTT into a risk-optimistic travel time index.

### 2.5 Combined mean travel time

As mentioned above, MBTT is a risk-optimistic travel time index reflecting the performance of reliable cases, whereas METT is just the opposite. Hence, MBTT and METT can be regarded as two boundaries (i.e. risk-optimistic and risk-pessimistic) of routing risk-attitudes. In reality, the extreme behaviours conducted by MBTT and METT may be not common. Most travellers consider both reliable aspect and unreliable aspect to make rational and overall decisions. In addition, due to heterogeneous personalities, travellers' routing risk-attitudes should be also various. Accordingly, we suggest using a combined mean travel time (CMTT) index to describe the heterogeneous routing risk-attitudes. Mathematically, the alpha-reliable CMTT here is formulated as the convex combination of alpha-reliable MBTT and METT, i.e.,

$$\psi_k^{rs}(\alpha;\lambda) = \lambda\varphi_k^{rs}(\alpha) + (1-\lambda)\eta_k^{rs}(\alpha) = \mu_k^{rs} + \frac{(\alpha-\lambda)\sigma_k^{rs}}{\sqrt{2\pi}\alpha(1-\alpha)}\exp\left(-\frac{\left[\Phi^{-1}(\alpha)\right]^2}{2}\right)\ \forall k \in K^{rs}, r \in R, s \in S \tag{4}$$

Here, the non-negative scalar $\lambda \in [0, 1]$ reflects the degree of risk-optimistic, and a larger value represents a more optimistic risk-attitude. According to Eq. (4), when $\lambda$ takes 0, $\alpha$, and 1, $\psi_k^{rs}(\alpha;\lambda)$ returns to the alpha-reliable METT, MTT, and MBTT, respectively. Therefore, CMTT is able to describe various routing risk-attitudes. The formulation of alpha-reliable METT $\eta_k^{rs}(\alpha)$ is supplied in Table 1.

### 2.6 Comparison analyses

In this sub-section, the routing risk-attitudes conducted by above mentioned four travel time indices (i.e.

MBTT, TTB/PTT, METT and CMTT) within a general framework will be compared. For this, a general alpha-reliable travel time index is defined, i.e.,

$$G_k^{rs}(\alpha) = \mu_k^{rs} + c\sigma_k^{rs} \ \forall k \in K^{rs}, r \in R, s \in S \tag{5}$$

Here, $c$ is an underdetermined coefficient, which reflects a traveller's risk-attitude towards variability. Risk-optimism, risk-neutralism and risk-pessimism are specified when $c$ takes negative value, zero, and positive value, respectively. Table 1 lists the formulations of $c$ for diverse travel time indices.

**Table 1** Routing risk-attitudes and formulations of $c$ for four travel time indices

| Alpha-reliable travel time index | Formulation of $c$ in Eq.(5) | Value of coefficient $c$ / Routing risk-attitude | | |
|---|---|---|---|---|
| | | $\alpha < 0.5$ | $\alpha = 0.5$ | $\alpha > 0.5$ |
| MBTT | $-\frac{1}{\sqrt{2\pi}\alpha}\exp\left(-\frac{\left[\Phi^{-1}(\alpha)\right]^2}{2}\right)$ | negative/ risk-optimistic | $-\sqrt{\frac{2}{\pi}}$ /risk-optimistic | negative/ risk-optimistic |
| PTT [1]/TTB [10] | $\Phi^{-1}(\alpha)$ | negative/ risk-optimistic | 0/ risk-neutral | positive/ risk-pessimistic |
| METT [17] | $\frac{1}{\sqrt{2\pi}(1-\alpha)}\exp\left(-\frac{\left[\Phi^{-1}(\alpha)\right]^2}{2}\right)$ | positive /risk-pessimistic | $\sqrt{\frac{2}{\pi}}$ /risk-pessimistic | positive/ risk-pessimistic |
| CMTT | $\frac{(\alpha-\lambda)}{\sqrt{2\pi}\alpha(1-\alpha)}\exp\left(-\frac{\left[\Phi^{-1}(\alpha)\right]^2}{2}\right)$ | negative/ risk-optimistic if $\lambda > \alpha$ , 0/ risk-neutral if $\lambda = \alpha$ , positive/ risk-pessimistic if $\lambda < \alpha$ | | |

Table 1 indicates that MBTT and METT are purely risk-optimistic and risk-pessimistic travel time indices, respectively. TTB seems to be able to express both risk-optimism and risk-pessimism. However, TTB needs to impose an unreliable reliability (e.g. smaller than 0.5) to formulate the risk-pessimism. This doesn't coincide with the reality since everyone usually preserves a relatively high reliability (e.g. larger than 0.8). Hence, TTB is actually a risk-pessimistic travel time index. By contrast, CMTT is capable of capturing various risk-attitudes within rational reliability levels. Therefore, CMTT is a complement and effective travel time index for modelling routing behaviours and implementing sensitivity analysis for risk-attitudes.

## 3 Combined mean traffic equilibrium model and solution algorithm

### 3.1 Combined mean traffic equilibrium model

Let $\boldsymbol{\delta} = \left[\delta_{ak}^{rs}\right]$ be the route-link incidence matrix. The feasible region $\boldsymbol{\Omega}_1$ described in terms of route flow can be described as follows.

$$\sum_k f_k^{rs} = q^{rs} \ \forall r \in R, s \in S \tag{6}$$

$$v_a = \sum_r \sum_s \sum_k f_k^{rs} \delta_{ak}^{rs} \ \forall a \in A \tag{7}$$

$$f_k^{rs} \ge 0 \ \forall k \in K^{rs}, r \in R, s \in S \tag{8}$$

Here, Eq. (6) is the conservation constraint of travel demand; Eq. (7) is a definitional conservation constraint between link flow and route flows; Eq. (8) promises route flow to be non-negative.

As is mentioned above, CMTT individuals aim to minimize the alpha-reliable CMTTs. These collective self-routing actions will finally reach a long-term habitual traffic equilibrium, which can be formulated as the following complementary expression.

$$f_k^{rs*}\left(\psi_k^{rs*}(\alpha)-\pi_\alpha^{rs}\right)=0,\ \psi_k^{rs*}(\alpha)-\pi_\alpha^{rs}\geq 0,\ f_k^{rs*}\geq 0, \forall k\in K^{rs}, r\in R, s\in S \tag{9}$$

The above expression is the alpha-reliable combined mean traffic equilibrium (CMTE) model. It implies $\psi_k^{rs*}(\alpha)=\pi_\alpha^{rs}$ if $f_k^{rs*}>0$ and $\psi_k^{rs*}(\alpha)\geq\pi_\alpha^{rs}$ if $f_k^{rs*}=0$, i.e., the used routes possess equal and minimum CMTT, while the unused ones possess the higher CMTTs. This equilibrium condition formulated in Eq. (9) is a Wardropian CMTE.

**Proposition 2.** The Wardropian CMTE expressed by Eq. (9) is equivalent to the following variational inequality problem (VIP).

$$\sum_{r\in R}\sum_{s\in S}\sum_{k\in K^{rs}}\psi_k^{rs*}(\alpha)\left(f_k^{rs}-f_k^{rs*}\right)\geq 0\ \forall \mathbf{f}\in\mathbf{\Omega_1} \tag{10}$$

with vector version expressed as

$$\boldsymbol{\psi}(\mathbf{f}^*)^{\mathrm{T}}\left(\mathbf{f}-\mathbf{f}^*\right)\geq 0\ \forall\mathbf{f}\in\mathbf{\Omega_1} \tag{11}$$

**Proof.** Proposition 2 can be guaranteed by the equivalent condition between nonlinear complementary problem (NCP) and VIP (see Facchinei and Pang [25]). □

**Proposition 3.** Assume CMTT is positive and continuous, and then the above CMTE model has at least one optimal solution.

**Proof.** According to Proposition 2, we only need to consider its equivalent VIP (11). Since $\mathbf{\Omega_1}$ is nonempty and convex, and the mapping $\boldsymbol{\psi}(\mathbf{f})$ is continuous. Hence, the VIP (11) has at least one solution (see Nagurney [26]). □

Consider the following link-route relationship

$$T_k^{rs}=\sum_{a\in A}T_a\delta_{ak}^{rs}\ \forall k\in K^{rs}, r\in R, s\in S \tag{12}$$

and the physical meaning of CMTT, it is reasonable to make the positive and continuous assumptions for CMTT in the above propositions. Therefore, the validity of VIP (11) and the existence of optimal solution are ensured.

### 3.2 Solution algorithm

Chen and Zhou [17] applied the modified alternating direction (MAD) algorithm [27] to solve the mean-excess traffic equilibrium (METE) model. The essence of MAD algorithm is transforming the linear constrained VIP into a uniform pattern, which helps to make the subsequent projection process much easier. Due to the economic meaning of Lagrangian multiplier vector $\boldsymbol{\pi}$ (i.e., the realised minimum cost vector between the O-D pairs) for the demand conservation constraints [25], VIP (11) can be reformulated as an equivalent VIP ($\mathbf{F}$, $\mathbf{\Omega}$) with closed and unified pattern as follows:

finding $\boldsymbol{u}^*\in\Omega$, such that

$$\mathbf{F}(\mathbf{u}^*)^{\mathrm{T}}(\mathbf{u}-\mathbf{u}^*)\geq\mathbf{0}\ \forall\mathbf{u}\in\mathbf{\Omega} \tag{13}$$

where

$$\mathbf{u}=\begin{pmatrix}\mathbf{f}\\ \boldsymbol{\pi}\end{pmatrix},\mathbf{F}(\mathbf{u})=\begin{pmatrix}\boldsymbol{\psi}(\mathbf{f})-\boldsymbol{\Lambda}^{\mathbf{T}}\boldsymbol{\pi}\\ \boldsymbol{\Lambda}^{\mathbf{T}}\mathbf{f}-\mathbf{Q}\end{pmatrix},\boldsymbol{\Omega}=\mathbf{R}_{+}^{\mathbf{m}}\times\mathbf{R}_{+}^{\mathbf{w}} \tag{14}$$

Here, $\boldsymbol{\Lambda}$ is the route-OD incidence matrix; $\mathbf{R}_{+}^{\mathbf{m}}$ and $\mathbf{R}_{+}^{\mathbf{w}}$ represent the *m*-dimensional and *w*-dimensional non-negative Euclidean spaces, respectively. $m$ and $w$ are the total numbers of routes and OD pairs, respectively. $\mathbf{Q}$ is a vector grouping travel demand $q^{rs}$.

The extra-gradient method [28], rather than the MAD algorithm [27], is employed here for two reasons. Firstly, the convergence condition for extra-gradient method is less strict [29]. Secondly, this study focuses on the illustration of the central idea of CMTT, and the extra-gradient method is adequate.

## 4 Numerical example

### 4.1 Network characteristics

The following example was conducted to perform the algorithm and characterize the CMTE model. For this reason, a hypothetical transportation network (see Fig. 1) with six paths from origin 1 to destination 10 was applied.

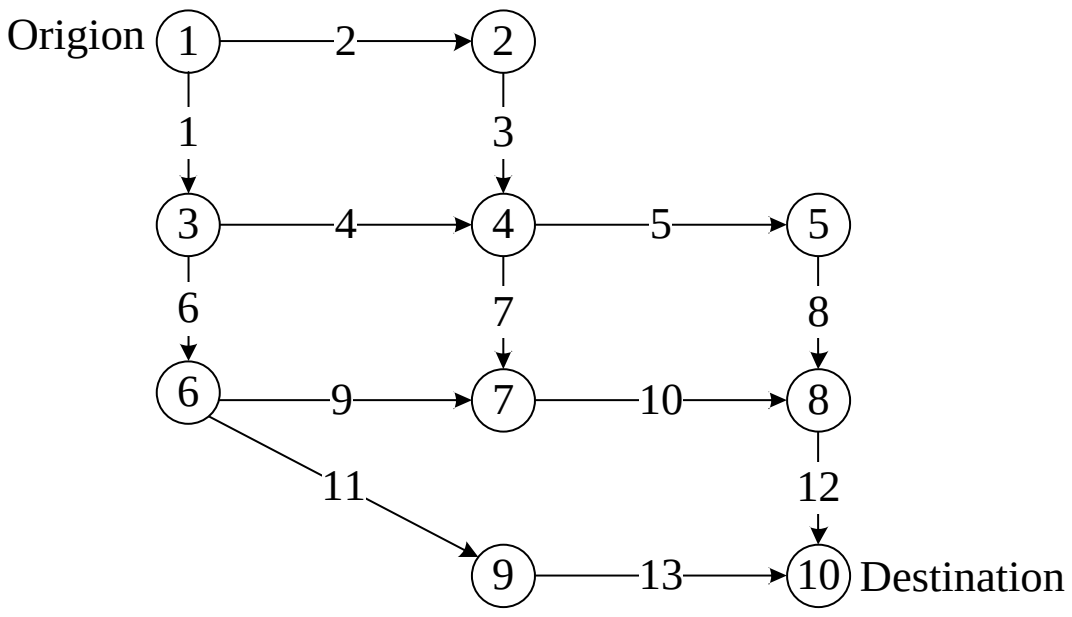

**Fig. 1** Test network

BPR link cost function was adopted in the subsequent numerical simulation.

$$T_a\left(v_a\right)=t_a^0\left[1+\beta\left(\frac{v_a}{C_a}\right)^n\right]\ \forall a\in A \tag{15}$$

where $t_a^0$ was the free-flow travel time on link $a$, $C_a$ was the link capacity, and $\beta$ and $n$ were the predefined parameters. Recently, based on the BPR formulation, many investigations were explored to study the stochastic traffic equilibrium under various disturbances, such as variation for link capacity [9−10] and travel demand [30]. Later, the variations considering both supply aspect and demand aspect were explored [15−16]. For simplicity and without loss of generality, only the link capacity degradation was considered here. Assumed the link capacities distributed uniformly, i.e., $C_a \sim U(\theta_a\overline{C}_a,\ \overline{C}_a)\ \forall a\in A$, where $\overline{C}_a$ and $\theta_a$ (vector $\boldsymbol{\Theta}=(\cdots,\theta_a,\cdots)$) were the design capacity and its degradable degree of link $a$, respectively. Then, the expectations and variations of link travel times could be deduced (see Lo and Tung [9] for detail).

$$E\left(T_a\right)=t_a^0+\beta t_a^0 v_a^n\frac{(1-\theta_a^{1-n})}{\overline{C}_a^n(1-\theta_a)(1-n)}\ \forall a\in A \tag{16}$$

$$Var\left(T_a\right)=\beta^2 (t_a^0)^2 v_a^{2n}\left\{\frac{(1-\theta_a^{1-2n})}{\overline{C}_a^{2n}(1-\theta_a)(1-2n)}-\left[\frac{(1-\theta_a^{1-n})}{\overline{C}_a^{n}(1-\theta_a)(1-n)}\right]^2\right\}\ \forall a\in A \tag{17}$$

where $\beta=0.15$ and $n=4$. The standard values of other parameters were listed in Table 2.

**Table 2** Network parameters and information

| Link | $t_a^0$ (min) | $\overline{C}_a$ (pcu/h) | Link | $t_a^0$ (min) | $\overline{C}_a$ (pcu/h) |
|---|---|---|---|---|---|
| 1 | 10 | 1000 | 8 | 10 | 1000 |
| 2 | 10 | 1000 | 9 | 4 | 1500 |
| 3 | 10 | 1000 | 10 | 10 | 2000 |
| 4 | 5 | 1600 | 11 | 30 | 1000 |
| 5 | 10 | 1000 | 12 | 10 | 1000 |
| 6 | 5 | 1000 | 13 | 10 | 1000 |
| 7 | 10 | 1000 | -- | -- | -- |

## 4.2 Numerical results

Based on Fig. 1, three experimental scenarios were designed.

**Scenario 1.** Convergence analysis of the algorithm.

**Scenario 2.** Sensitivity analysis (SA) for risk-attitudes under different demand levels.

**Scenario 3.** SA for the risk-attitudes under different capacity degradation levels.

The numerical parametric values of the above three scenarios were listed in Table 3, where $Q$ denoted the travel demand and vector $\mathbf{\Theta}$ denoted the capacity degradable degree. For simplicity, $\mathbf{\Theta}=0.8$ meant $\theta_a=0.8\ \forall a\in A$.

**Table 3** Numerical scenarios

| Scenario | $\alpha$ | $Q$ (pcu/h) | $\mathbf{\Theta}$ | $\lambda$ |
|---|---|---|---|---|
| 1 | 0.9 | 4000 | 0.8 | 0.5 |
| 2 | 0.9 | [3000:1000:4000] | 0.8 | [0.0:0.1:1.0] |
| 3 | 0.9 | 4000 | [0.6:0.1: 0.9] | [0.0:0.1:1.0] |

The numerical results for Scenario 1, 2 and 3 are displayed in Fig. 2, 3 and 4, respectively. In Fig. 2 and 3, the indicator $\lambda_0$ denotes the travelers across the network are risk-neutral at point $\lambda=\lambda_0$. Fig. 2 performs the evolution of the average network travel time (ANTT). Note that, in the following context, all units will be omitted to simplify the presentation.

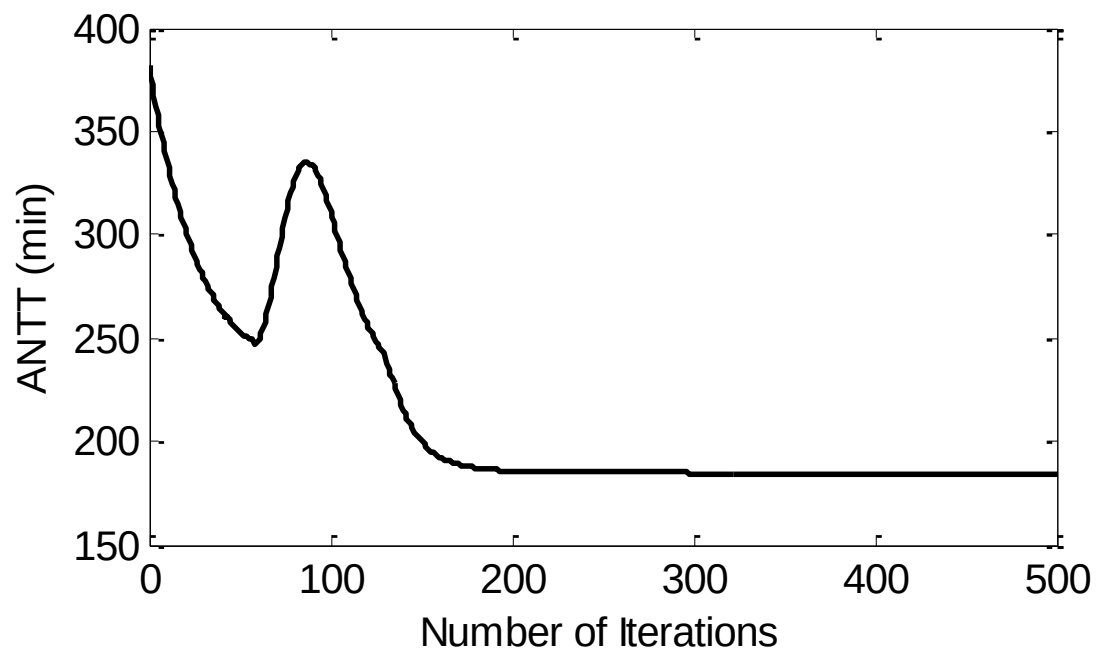

**Fig. 2** Convergence curve

Fig. 2 shows ANTT fluctuates obviously in the beginning 100 iterations, and decreases steadily until it converges to the required accuracy (equalling to 1.0e-4) as iteration proceeds. The whole computation process, executed on a microcomputer with CPU frequency of 1.5G and 512M RAM, consumes about 2 seconds. Consequently, the current solution algorithm is adequate to the following simulations.

Fig. 3 displays the effects of travellers' risk-attitudes on the network performance (measured by ANTT) under different demand levels.

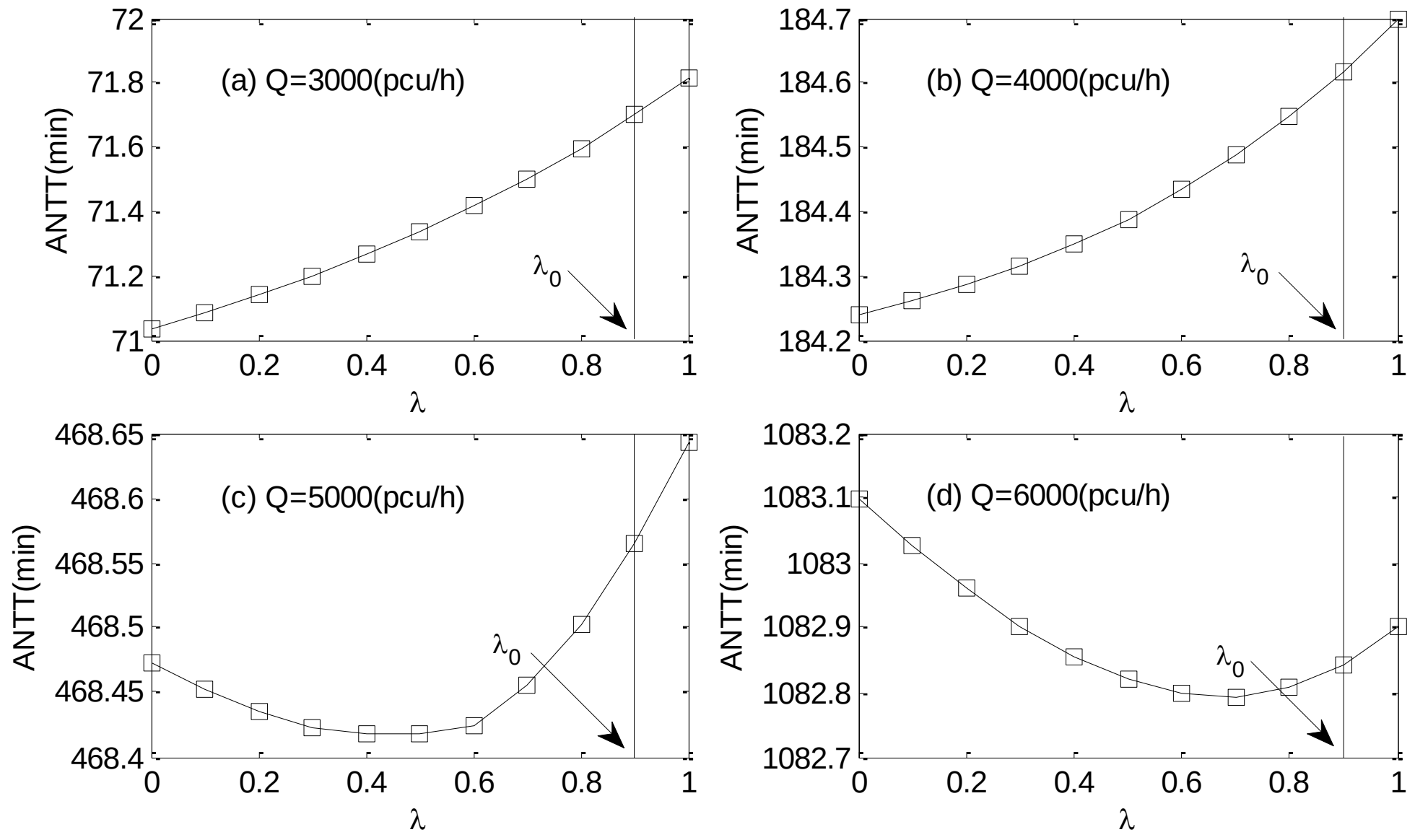

**Fig. 3** SA for travelers' risk-attitudes under different demand levels

Fig. 3 indicates i) ANTT surges to about 1083 from about 72 as $Q$ increases from 3000 to 6000; ii) under a given demand level, the variation of ANTT is not obvious as $\lambda$ increases, which provides some evidence for the existence of risk-optimistic routing behaviours, because risk-optimists do not feel apparent loss from their optimistic routing behaviour in the long run; iii) Overall, modest risk-pessimist helps to alleviate the network congestion and save personal MTT. In addition, when the demand level is low (e.g. $Q \le 4000$), METT travellers enjoy less personal MTT than MBTT ones, whereas it is just the opposite for high demand level cases (e.g. $Q \ge 5000$). This demonstrates that the optimists do not always suffer from the most loss in daily travels, and this phenomenon is especially obvious in the high demand cases.

Fig. 4 displays the effects of travellers' risk-attitudes on the network performance (measured by ANTT) under different capacity degradation levels.

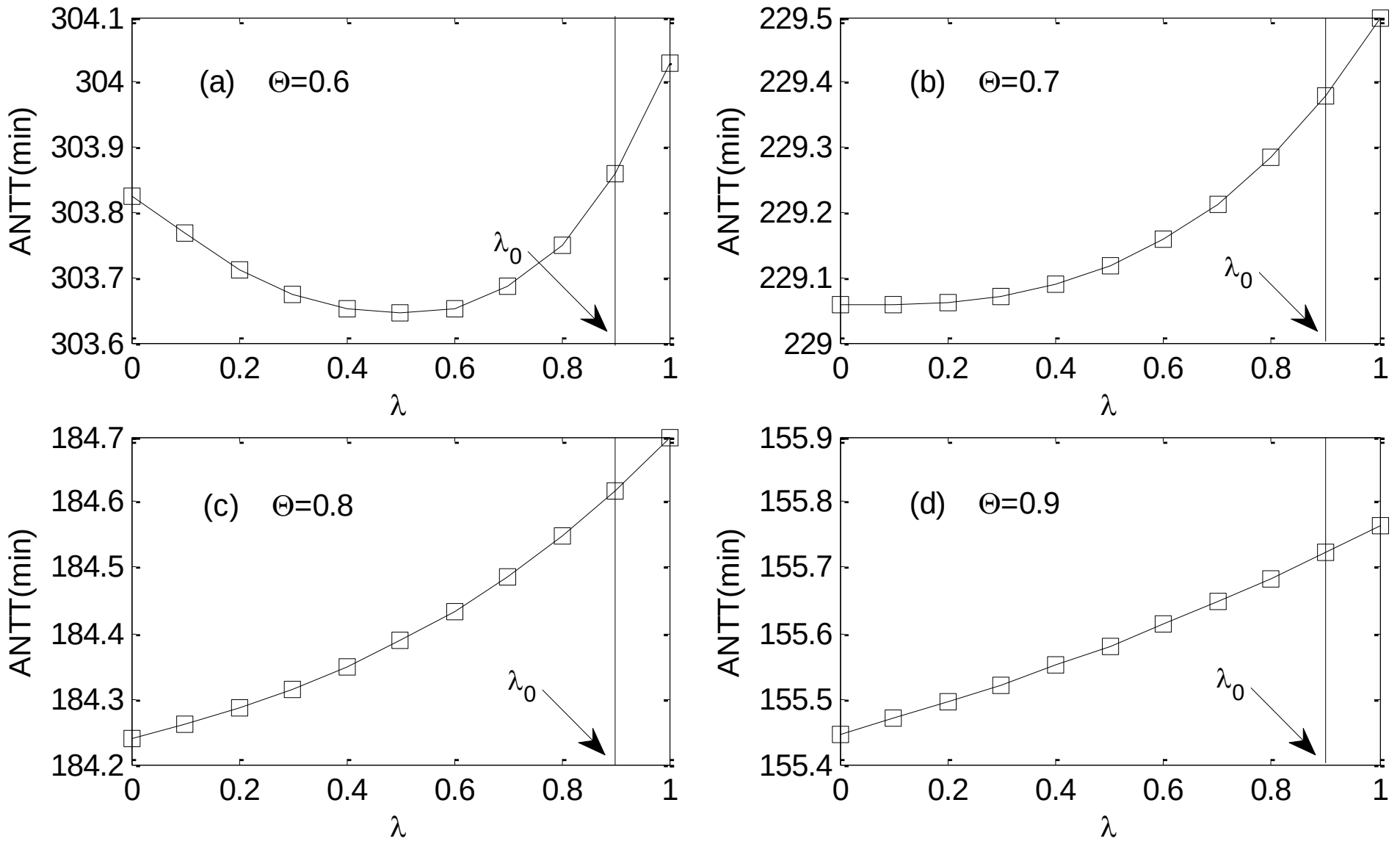

**Fig. 4** SA for travelers' risk-attitudes under different degradable levels of capacity

Fig. 4 shows some analogous phenomena as Fig. 3 does, such as i) ANTT decreases to about 156 from about 304 as $\boldsymbol{\Theta}$ increases from 0.6 to 0.9; ii) under a specific capacity degradation level, unobvious changes happen on ANTT as $\lambda$ increases, which further verifies the existence of risk-optimistic routing behaviours; iii) Generally, Modest risk-pessimist contributes to alleviating the network congestion and thrift of private MTT. Moreover, when the capacity degradation is relatively low (e.g. $\boldsymbol{\Theta} \geq 0.8$), METT travellers enjoy less private MTT than MBTT travellers, while the superiority decreases in the higher degradation level cases (e.g. $\boldsymbol{\Theta} \leq 0.7$).

## 5 Conclusions and future research

1) Compare with the existing travel time indices (e.g. PTT/TTB, METT), CMTT is capable of capturing various routing risk-attitudes, spreading from risk-optimism to risk-pessimism.

2) Numerical studies demonstrate the efficiency of algorithm; it can be also founded that, for low and modest demand and capacity degradation levels, the risk-pessimistic routing behaviour is more beneficial to both network and individuals; however, as the demand and capacity degradation levels increase, the superiority of risk-pessimism becomes more and more obvious.

Many further works are worthy of exploring based on the proposed CMTE model.

1) More efficient algorithms for generating the CMTT routes and solving the proposed model need to be developed and tested on a large-scale network.

2) Empirical studies need to be performed to obtain a better understanding of the travellers' routing risk-attitudes under the traffic information service (e.g. Liu et al. [31]).

3) To incorporate multiple risk-attitudes into the current CMTE model is also a valuable extension;

4) Finally, by factoring various sources of travel time variation and general cost into the proposed modelling approach, the CMTE model will be more realistic.

## Appendix A: Estimating MBTT

Note that some super-scripts and sub-scripts will be omitted to facilitate the following estimation. Since $T \sim N(\mu,\sigma)$ and $\alpha = \Pr\left(T \leq \xi(\alpha)\right)$, the moment generating function $M_{T:T\leq\xi(\alpha)}(z)$ can be formulated as

follows:

$$M_{T:T\le\xi(\alpha)}(z)=\int_{-\infty}^{\xi(\alpha)}e^{tz}\frac{1}{\Pr\left(T\le\xi(\alpha)\right)}\frac{1}{\sqrt{2\pi}\sigma}e^{-\frac{(t-\mu)^2}{2\sigma^2}}dt. \tag{A.1}$$

For brevity, we rewrite $M_{T:T\le\xi(\alpha)}(z)$ as $M(z)$ in the remaining deduction. Let $T=\mu+\sigma S$, where $S\sim N(0,1)$, then Eq. (A.1) can be further deduced as

$$\begin{aligned}M(z)&=\frac{1}{\alpha}e^{\mu z}\int_{-\infty}^{\frac{\xi(\alpha)-\mu}{\sigma}}e^{\sigma sz}\frac{1}{\sqrt{2\pi}}e^{-\frac{s^2}{2}}ds\\&=\frac{1}{\alpha}e^{\mu z}\int_{-\infty}^{\frac{\xi(\alpha)-\mu}{\sigma}}e^{\frac{\sigma^2z^2}{2}}\frac{1}{\sqrt{2\pi}}e^{-\frac{(s-\sigma z)^2}{2}}ds\\&=\frac{1}{\alpha}\exp\left(\mu z+\frac{\sigma^2z^2}{2}\right)\int_{-\infty}^{\frac{\xi(\alpha)-\mu}{\sigma}}\frac{1}{\sqrt{2\pi}}e^{-\frac{(s-\sigma z)^2}{2}}ds.\end{aligned}$$

Let $Y=S-\sigma z$, then

$$M(z)=\frac{1}{\alpha}\exp\left(\mu z+\frac{\sigma^2z^2}{2}\right)\int_{-\infty}^{\frac{\xi(\alpha)-\mu}{\sigma}-\sigma z}\frac{1}{\sqrt{2\pi}}e^{-\frac{y^2}{2}}dy=\frac{1}{\alpha}\exp\left(\mu z+\frac{\sigma^2z^2}{2}\right)\Phi\left(\frac{\xi(\alpha)-\mu}{\sigma}-\sigma z\right). \tag{A.2}$$

Taking the first-order derivative of Eq. (A.2) with respect to *z*, we obtain

$$M'(z)=\frac{1}{\alpha}(\mu+\sigma^2z)\exp\left(\mu z+\frac{\sigma^2z^2}{2}\right)\Phi\left(\frac{\xi(\alpha)-\mu}{\sigma}-\sigma z\right)-\frac{\sigma}{\alpha}\exp\left(\mu z+\frac{\sigma^2z^2}{2}\right)\phi\left(\frac{\xi(\alpha)-\mu}{\sigma}-\sigma z\right), \tag{A.3}$$

where $\Phi(y)=\int_{-\infty}^{y}\frac{1}{\sqrt{2\pi}}e^{-\frac{x^2}{2}}dx,\ \phi(y)=\frac{1}{\sqrt{2\pi}}e^{-\frac{y^2}{2}}$. Let z=0 in Eq. (A.3) and define $T=\mu+\sigma S$, because $\alpha\equiv\int_{-\infty}^{\xi(\alpha)}f(t)dt$, we have $\alpha\equiv\int_{-\infty}^{\frac{\xi(\alpha)-\mu}{\sigma}}\phi(s)ds\equiv\Phi\left(\frac{\xi(\alpha)-\mu}{\sigma}\right)$ and

$$M'(0)=\mu-\frac{\sigma}{\alpha}\phi\left(\frac{\xi-\mu}{\sigma}\right). \tag{A.4}$$

Let $\Phi^{-1}(\alpha)=\frac{\xi-\mu}{\sigma}$, we finally derive the alpha-reliable mean-below travel time (MBTT) index as follows.

$$\varphi(\alpha)=M'(0)=\mu-\frac{\sigma}{\sqrt{2\pi}\alpha}\exp\left(\frac{\left[\Phi^{-1}(\alpha)\right]^2}{2}\right). \tag{A.5}$$

## Reference

[1] FRANK H. Shortest paths in probabilistic graphs [J]. Operations Research, 1969, 17(4): 583–599.

[2] FAN Y, KALABA R, MOORE J. Arriving on time [J]. Journal of Optimization Theory and Application, 2005, 127(3): 497–513.

[3] LEVY H. Stochastic dominance: Investment decision making under uncertainty [M]. New York: Springer, 2006: 147–152.

[4] NIE Y, WU X. Shortest path problem considering on-time arrival probability [J]. Transportation Research Part B, 2009, 43(6): 597–613.

[5] NIE Y, WU X. Reliable a priori shortest path problem with limited spatial and temporal dependencies [C]// Proceedings of the 18th International Symposium on Transportation and Traffic Theory, Hong Kong, 2009: 169–195.

[6] NIE Y. Multi-class percentile user equilibrium with flow-dependent stochasticity [J]. Transportation Research Part B, 2011, 45(10): 1641–1659.

[7] HALL R W. The fastest path through a network with random time-dependent travel time [J]. Transportation Science,

1986, 20(3): 182–186.

[8] WU X, NIE Y. Modelling heterogeneous risk-taking behaviour in route choice: A stochastic dominance approach [J]. Transportation Research A, 2011, 45(9): 896–915.

[9] LO H K, TUNG Y K. Network with degradable links: Capacity analysis and design [J]. Transportation Research Part B, 2003, 37(4): 345–363.

[10] LO H K, LUO X W, SIU B W Y. Degradable transport network: Travel time budget of travellers with heterogeneous risk aversion [J]. Transportation Research Part B, 2006, 40(9): 792–806.

[11] ARTZNER P, DELBAEN F, EBER J M, Heath D. Coherent measures of risk [J]. Mathematical Finance, 1999, 9(3): 203–228.

[12] CHEN A, JI Z W. Path finding under uncertainty [J]. Journal of Advanced Transportation, 2005, 39(1): 19–37.

[13] SHAO H, LAM W H K, TAM M L. A reliability-based stochastic traffic assignment model for network with multiple user classes under uncertain in demand [J]. Network and Spatial Economics, 2006, 6(3/4): 313–332.

[14] SIU B W Y, LO H K. Doubly uncertain transport network: Degradable link capacity and perception variation in traffic conditions [J]. Transportation Research Record, 2006, 1964: 59–69.

[15] SHAO H, LAM W H K, TAM M L, YUAN X M. Modelling rain effects on risk-taking behaviours of multi-user classes in road networks with uncertainty [J]. Journal of Advanced Transportation, 2008, 42(3): 265–290.

[16] LAM W H K, SHAO H, SUMALEE A. Modelling impacts of adverse weather conditions on a road network with uncertainties in demand and supply [J]. Transportation Research Part B, 2008, 42(10): 890–910.

[17] CHEN A, ZHOU Z. The alpha-reliable mean-excess traffic equilibrium model with stochastic travel times [J]. Transportation Research Part B, 2010, 44(4): 493–513.

[18] ROCKAFELLAR R T, URYASEV S. Optimization of conditional value-at-risk [J]. Journal of Risk, 2000, 2(3): 21–41.

[19] FOSGERAU M, KARLSTROM A. The value of reliability [J]. Transportation Research Part B, 2010, 44(1): 38–49.

[20] FOSGERAU M, ENGELSON L. The value of travel time variance [J]. Transportation Research Part B, 2011, 45(1): 1–8.

[21] ORDONEZ F, STIER-MOSES N E. Wardrop equilibria with risk-averse users [J]. Transportation Science, 2010, 44(1): 63–86.

[22] NG M W, SZETO W Y, WALLER S T. Distribution-free travel time reliability assessment with probability inequalities [J]. Transportation Research Part B, 2011, 45(6): 852–866.

[23] GERVAIS S, HEATON J B, ODEAN T. Overconfidence, Investment Policy, and Executive Stock Options [R]. Rodney L. White Center for Financial Research Working Paper No. 15-02, 2002.

[24] BILINGSLEY P. Probability and measure [M]. New York: John Wiley and Sons, Inc., 1995: 359–363.

[25] FACCHINEI F, PANG J S. Finite-dimensional variational inequalities and complementarity problems [M]. New York: Springer, 2003: 41–46.

[26] NAGURNEY A. Network economics: A variational inequality approach [M]. Dordrecht: Kluwer Academic Publishers, 1993: 14–15.

[27] HAN D. A modified alternating direction method for variational inequation problems [J]. Applied Mathematics and Optimization, 2002, 45(1): 63–74.

[28] KORPELEVICH G M. The extra-gradient method for finding saddle points and other problems [J]. Matecon, 1976, 12: 747–756.

[29] HAN Ji-ye, XIU Nai-hua, QI Hou-duo. Nonlinear complementarity theory and algorithms [M]. Shanghai: Shanghai scientific & Technical Publishers, 2006: 132–141. (in Chinese)

[30] ZHOU Z, CHEN A. Comparative analysis of three user equilibrium models under stochastic demand [J]. Journal of Advanced Transportation, 2008, 42(3): 239–263.

[31] LIU Tian-liang, HUANG Hai-jun, TIAN Li-jun. Microscopic simulation of multi-lane traffic under dynamic tolling and information feedback [J]. Journal of Central South University of Technology, 2009, 16(5): 0865–0870.